\documentclass[letterpaper, 10pt, conference]{ieeeconf}
 \IEEEoverridecommandlockouts 
\usepackage{setspace, epsfig, psfrag, amssymb, amsmath, cite} \newfont{\smalll}{cmr8}

\def\IR{\mathbb{R}}
\def\IS{\mathbb{S}}
\def\zero{\mathbf{0}}
\def\one{\mathbf{1}}
\def\spec{\mathrm{spec}}
\def\diag{\mathrm{diag}}
\def\span{\mathrm{span}}
\def\eye{\mathrm{I}}

\def\IS{\hbox{I\hskip-.1em S}}
\def\IC{\hbox{C\hskip-
.5em\raise.5ex\hbox{$\scriptscriptstyle\mid$}}\ }
\def\Ic{\hbox{\smalll C\hskip-
.5em\raise.3ex\hbox{$\scriptscriptstyle\mid$}}\ }

\def\T={\buildrel {\scriptscriptstyle\triangle} \over =}

\def\sqr#1#2{{\vcenter{\vbox{\hrule height.#2pt\hbox{\vrule
width.#2pt height#1pt \kern#1pt\vrule width.#2pt}\hrule
height.#2pt}}}}
\def\square{\mathchoice\sqr64\sqr64\sqr33\sqr33}
\def\diag{\mathop{\rm diag}}

\def\block-diag{\mathop{\rm block{\scriptstyle -}diag}}

\def\pmbb#1{\setbox0=\hbox{#1}\raise 0.5ex\box0}

\def\norm#1{\|#1\|}

\newtheorem{lemma}{Lemma}[section]
\newtheorem{theorem}{Theorem}[section]

\newtheorem{remark}{Remark}[section]

\newtheorem{proposition}{Proposition}[section]

\newtheorem{assumption}{Assumption}[section]

\newcommand{\bequ}{\begin{eqnarray}}
\newcommand{\eequ}{\end{eqnarray}}
\newcommand{\mT}{^\mathrm{T}}

\newcommand{\rom}{\mathrm}

\newcommand {\teq}      {\triangleq}

\newcommand {\beq}      {\begin{equation}}
\newcommand {\eeq}      {\end{equation}}

\def\IR{{\mathbb R}}
\def\IC{{\mathbb C}}

\def\IS{{\mathbb S}}

\renewcommand\theenumi{\roman{enumi}}
\renewcommand\theenumii{\alph{enumii}}
\renewcommand\theenumiii{\arabic{enumiii}}
\renewcommand\theenumiv{\Alph{enumiv}}

\def\sc{s_{\rm c}}

\def\sc{s_{\rm c}}

\begin{document}
    {
    \title{\Huge{ {Distributed Control of Multiagent Systems\\under Unknown Persistent Disturbances}}}
    }{
    \author{{Tansel Yucelen and Magnus Egerstedt 
    }
    \thanks{This is the extended version of the paper ``Control of Multiagent Systems under Persistent Disturbances'' published at the American Control Conference 2012 by T. Yucelen and M. Egerstedt.}
    \thanks{T. Yucelen is with the Department of Mechanical and Aerospace Engineering, Missouri University of Science and Technology, Rolla, MO 65401, E-mail: {\tt\small yucelen@mst.edu}.}
    \thanks{M. Egerstedt is with the School of Electrical and Computer Engineering, Georgia Institute of Technology, Atlanta, GA 30332, E-mail: {\tt\small magnus@ece.gatech.edu}.}
    }
    }
\markboth{} {Shell \MakeLowercase{\textit{et al.}}: Bare Demo of IEEEtran.cls for Journals} \newcommand{\eqnref}[1]{(\ref{#1})}
\newcommand{\class}[1]{\texttt{#1}} \newcommand{\package}[1]{\texttt{#1}} \newcommand{\file}[1]{\texttt{#1}} \newcommand{\BibTeX}{\textsc{Bib}\TeX}
\maketitle \baselineskip 12.5pt \parskip 4pt \maketitle
\begin{abstract}
This paper focuses on the consensus and formation problems of multiagent systems under unknown persistent disturbances. 
Specifically, we propose a novel method that combines an existing consensus (or formation) algorithm with a new controller. 
The new controller has an integral action that produces a control input based on an error signal locally projected onto the column space of the graph Laplacian. 
This action allows agents to achieve a consensus or a predetermined formation objective under constant or time-varying disturbances. 
This study has analyzed the stability properties of this architecture and provides several illustrative examples to demonstrate the proposed approach.

\end{abstract}
\vspace{-0.1cm}
\section{Introduction}
\vspace{-0.1cm}

Numerous studies from diverse areas of science and engineering have conducted research in the recent decades that focused on the consensus and formation problems in multiagent systems
(see \cite{ref:002,ref:a01,ref:a02} and references therein).
Specifically, consensus refers to agents coming to a global agreement on a state value, and formation refers to agents moving toward a desired geometric shape.
In this paper, we explore new solutions to these problems for the case when the dynamics of each agent is perturbed by an unknown persistent disturbance.

We propose an architecture that combines an existing consensus (or formation) algorithm with a new controller.
This controller's integral action produces a control input based on an error signal that is locally projected onto the column space of the graph Laplacian,
and this action allows agents to achieve a consensus or a predetermined formation objective under constant or time-varying disturbances.
In particular, local projection removes the effect of the bias that could possibly contained in this error signal,
and hence, this controller does not necessarily drive the states of the agents to the origin while suppressing the effect of disturbances.
The projected error signal is constructed by taking the difference between an agent's state and its predicted state.
In particular, the predicted state is obtained from a differential equation that resembles an ideal consensus (or formation) model modified by this error term in such a way that
the ideal consensus (or formation) model is recovered (resp., approximately recovered) when this error term goes to zero (resp., goes to a small neighborhood around zero).
We refer to \cite{ref:005} for more details regarding the state predictor which is used in the context of adaptive control theory.
The realization of the proposed controller only requires an agent to have access to its own state and to relative state information with respect to its neighbors.
We verify the stability properties of this architecture by using the results from linear algebra, matrix mathematics, and the Lyapunov theory.

The organization of the paper is as follows.
Section II recalls some of the basic notions from graph theory and then provides the general setup of the consensus and formation problems for multiagent networks.
Section III presents the formulation for the consensus problem.
Section IV provides the stability analysis for the architecture proposed in Section III, including extensions to the formation problem.
Section V modifies the proposed approach of Section III for cases in which the agents need to come to a global agreement on a constant point in space.
Section VI extends the results of all previous sections to include the time-varying disturbance case.
Several illustrative examples on a cycle graph with six agents are given in Section VII to demonstrate the efficacy of the proposed approach,
and conclusions are summarized in Section VIII.

The notation used in this paper is fairly standard.
Specifically, $\IR$ denotes the set of real numbers,
$\IR^n$ denotes the set of $n \times 1$ real column vectors,
$\IR^{n \times m}$ denotes the set of $n \times m$ real matrices,
$\IR_+$ denotes the set of positive real numbers,
$\IR_+^{n \times n}$ (resp., $\overline{\IR}_+^{\hspace{0.1em} n \times n}$) denotes the set of $n \times n$ positive-definite (resp., nonnegative-definite) real matrices,
$\IS_+^{n \times n}$ (resp., $\overline{\IS}_+^{\hspace{0.1em} n \times n}$) denotes the set of $n \times n$ symmetric positive-definite (resp., symmetric nonnegative-definite) real matrices,
$\zero_n$ denotes the $n \times 1$ vector of all zeros,
$\one_n$ denotes the $n \times 1$ vector of all ones,
$\zero_{n \times n}$ denotes the $n \times n$ zero matrix,
and $\eye_n$ denotes the $n \times n$ identity matrix.
Furthermore, we write $(\cdot)\mT$ for transpose,
$(\cdot)^{-1}$ for inverse,
$(\cdot)^+$ for the Moore-Penrose generalized inverse,
$\otimes$ for the Kronecker product,
$\norm{\cdot}_2$ for the Euclidian norm,
$\norm{\cdot}_\rom{F}$ for the Frobenius matrix norm,
$\lambda_\rom{min}(A)$ (resp., $\lambda_\rom{max}(A)$) for the minimum (resp., maximum) eigenvalue of the Hermitian matrix $A$,
$\lambda_i(A)$ for the $i$-th eigenvalue of $A$ ($A$ is symmetric and the eigenvalues are ordered from least to greatest value),
$\spec(A)$ for the spectrum of $A$,
$\pi_0(A)$, $\pi_+(A)$, and $\pi_-(A)$ for the number of eigenvalues (counted with algebraic multiplicities) of $A$ having zero, positive, and negative real parts, respectively,
$\diag(a)$ for the diagonal matrix with the vector $a$ on its diagonal,
$[A]_{ij}$ for the entry of the matrix $A$ on the $i$-th row and $j$-th column,
$\mathcal{N}(A)$ for the null space of the matrix $A$,
and $\span(a)$ for the span of the vector $a$ (the subspace generated by scalar multiplies of $a$).

\vspace{-0.1cm}
\section{Preliminaries}
In this section we recall some of the basic notions from graph theory, which is followed by the general setup of the consensus and formation problems for multiagent networks.
We refer to \cite{ref:001} and \cite{ref:002} for more details about graph theory and multiagent networks.

\vspace{-0.1cm}

\subsection{Graphs and Their Algebraic Representation}
In the multiagent literature, graphs are broadly adopted to encode interactions in networked systems.
An \textit{undirected} graph $\mathcal{G}$ is defined by a set $\mathcal{V}_\mathcal{G}=\{1,\ldots,n\}$ of \textit{nodes}
and a set $\mathcal{E}_\mathcal{G} \subset \mathcal{V}_\mathcal{G} \times \mathcal{V}_\mathcal{G}$ of \textit{edges}.
If $(i,j) \in \mathcal{E}_\mathcal{G}$, then the nodes $i$ and $j$ are \textit{neighbors} and the neighboring relation is indicated with $i \sim j$.
The \textit{degree} of a node is given by the number of its neighbors.
Letting $d_i$ be the degree of node $i$, then the \textit{degree} matrix of a graph $\mathcal{G}$, $\mathcal{D}(\mathcal{G}) \in \IR^{n \times n}$, is given by
\bequ
        \mathcal{D}(\mathcal{G}) \triangleq \diag(d), \quad d=[d_1,\ldots,d_n]\mT. \label{DegMat}
\eequ
A \textit{path} $i_0 i_1 \ldots i_L$ is a finite sequence of nodes such that $i_{k-1} \sim i_k$, $k=1, \ldots, L$,
and a graph $\mathcal{G}$ is \textit{connected} if there is a path between any pair of distinct nodes.
The \textit{adjacency} matrix of a graph $\mathcal{G}$, $\mathcal{A}(\mathcal{G}) \in \IR^{n \times n}$, is given by
\begin{eqnarray}
   [\mathcal{A}(\mathcal{G})]_{ij}
   \teq
   \left\{ \begin{array}{cl}
      1, & \mbox{ if $(i,j)\in\mathcal{E}_\mathcal{G}$},
      \\
      0, & \mbox{otherwise}.
   \end{array} \right.
   \label{AdjMat}
\end{eqnarray}
The \textit{Laplacian} matrix of a graph, $\mathcal{L}(\mathcal{G}) \in \overline{\IS}_+^{\hspace{0.1em} n \times n}$, playing a central role in many graph theoretic treatments of multiagent systems is given by
\bequ
    \mathcal{L}(\mathcal{G}) \triangleq \mathcal{D}(\mathcal{G}) - \mathcal{A}(\mathcal{G}), \label{Laplacian}
\eequ
where the spectrum of the Laplacian for a connected, undirected graph can be ordered as
\bequ
    0 = \lambda_1(\mathcal{L}(\mathcal{G}))<\lambda_2(\mathcal{L}(\mathcal{G}))\le \cdots \le \lambda_n(\mathcal{L}(\mathcal{G})), \label{LapSpec}
\eequ
with $\one_n$ as the eigenvector corresponding to the zero eigenvalue $\lambda_1(\mathcal{L}(\mathcal{G}))$
and $\mathcal{L}(\mathcal{G}) \one_n = \zero_n$ and $\rom{e}^{\mathcal{L}(\mathcal{G})} \one_n = \one_n$ hold.

\subsection{Consensus Dynamics}
We can model a given multiagent system by a graph $\mathcal{G}$ where nodes and edges represent agents and interagent information exchange links,
respectively\footnote{Throughout this paper we assume that the network is static. As such, the movements of the agents will not cause edges to appear or disappear in the network}.
Let $x_i(t) \in \IR^m$, denote the state of node $i$ at time $t \ge 0$, whose dynamics is described by the single integrator
\bequ
        \dot{x}_i(t) = u_i(t), \quad x_i(0)=x_{i0}, \quad i=1,\cdots,n, \label{eqn:01}
\eequ
with $u_i(t) \in \IR^m$ being the control input of node $i$.
Allowing agent $i$ to have access to the relative state information with respect to its neighbors, the solution of the consensus problem can be achieved by applying
\bequ
            u_i(t)=-\sum_{i \sim j} \bigl(x_i(t)-x_j(t) \bigl), \label{eqn:02}
\eequ
to the single integrator dynamics given by (\ref{eqn:01}), where (\ref{eqn:01}) in conjunction with (\ref{eqn:02}) can be represented as the Laplacian dynamics of the form
\bequ
        \dot{x}(t)=-\mathcal{L}(\mathcal{G})\otimes \eye_m \hspace{0.2em} x(t), \quad x(0)=x_0, \label{eqn:03}
\eequ
where $x(t)=[x_1\mT(t),\cdots,x_n\mT(t)]\mT$ denotes the aggregated state vector of the multiagent system.
Although our results can directly be extended to the case of (\ref{eqn:03}), in what follows we will focus on the system
\bequ
        \dot{x}(t)=-\mathcal{L}(\mathcal{G})x(t), \quad x(0)=x_0, \label{eqn:04}
\eequ
capturing the multiagent dynamics with individual agent states evolving in $\IR$.
Furthermore, we assume that $\mathcal{G}$ is connected.
Considering (\ref{eqn:04}), we note that $x(t) \rightarrow \one_n (\one_n\mT \one_n)^{-1} \one_n\mT x_0 = ({\one_n \one_n\mT}/n)x_0$ as $t \rightarrow \infty$,
since the (undirected) graph $\mathcal{G}$ is connected.

\subsection{Formation Dynamics}
For our take on the formation problem, define $\tau_i$ as the displacement of $x_i$ from the target location $\zeta_i$.
Then, using the state transformation given by
\bequ
        \tau_i(t)=x_i(t)-\zeta_i, \quad i=1,\ldots,n, \label{eqn:05}
\eequ
the solution of the \textit{invariant} formation problem follows from (\ref{eqn:04}) as
\bequ
        \dot{x}(t)=-\mathcal{L}(\mathcal{G})x(t)-\mathcal{L}(\mathcal{G})\zeta, \quad x(0)=x_0, \label{eqn:06}
\eequ
where $\zeta=[\zeta_1, \cdots, \zeta_n]\mT$.
Note that (\ref{eqn:06}) can equivalently be written as
 \bequ
        \dot{x}_i(t) = -\sum_{i \sim j} \bigl(x_i(t)-x_j(t) \bigl)-\bigl(\zeta_i-\zeta_j \bigl), \ x_i(0)=x_{i0}. \label{eqn:07}
\eequ

\subsection{Lyapunov and Asymptotic Stability}
We conclude this section by noting that a matrix $A\in\IR^{n \times n}$ is Lyapunov stable if and only if every eigenvalue of $A$ has nonpositive real part
and every eigenvalue of $A$ with zero real part is semisimple \cite{ref:003},
and it is asymptotically stable (Hurwitz) if and only if every eigenvalue of $A$ has negative real part \cite{ref:004}.
Furthermore, the zero solution $x(t)\equiv0$ to a linear dynamical system of the form $\dot{x}(t)=Ax(t)$, $x(0)=x_0$, $t \ge 0$, is Lyapunov stable
if and only if $A$ is Lyapunov stable, and asymptotically stable if and only if $A$ is asymptotically stable (see, for example, Theorem 3.15 of \cite{ref:004}).
For a connected, undirected graph, $-\mathcal{L}(\mathcal{G})$ is Lyapunov stable, which is a direct consequence of (\ref{LapSpec}).
Therefore, the zero solution $x(t)\equiv0$ to (\ref{eqn:04}) is Lyapunov stable.

\section{Problem Formulation}
Consider a system of $n$ agents exchanging information among each other using their local measurements, according to a connected, undirected graph $\mathcal{G}$.
Let $x(t)=[x_1(t),\ldots,x_n(t)]\mT \in \IR^n$ and $u(t)=[u_1(t),\ldots,u_n(t)]\mT \in \IR^n$ denote the aggregated state vector and control input of the multiagent system at time $t \ge 0$, respectively.
Suppose that the dynamics of agents are perturbed by an \textit{unknown} \textit{constant} disturbance vector $w=[w_1,\ldots,w_n]\mT \in \IR^n$,
where $w_i$, $i =1, \ldots, n$, represents a disturbance effecting the dynamics of $i$-th agent\footnote{For the case when there is no disturbance effect on the dynamics of an agent $i$, $w_i=0$.
Note that we present extensions to time-varying disturbance case later.}.
Then, the multiagent system is represented by
\bequ
    \dot{x}(t)=u(t)+w, \quad x(0)=x_0. \label{pf:01}
\eequ

\subsection{Proposed Controller}
To address the consensus problem for the multiagent system given by (\ref{pf:01}), we propose the controller of the form
\bequ
    u(t) = u_\rom{s}(t) + u_\rom{a}(t), \label{pf:02}
\eequ
where $u_\rom{s}(t) \in \IR^n$ denotes the standard control input given by
\bequ
    u_\rom{s}(t) = -\mathcal{L}(\mathcal{G})x(t), \label{pf:03}
\eequ
and $u_\rom{a}(t) \in \IR^n$ denotes the additional control input given by
\bequ
    u_\rom{a}(t)= -\hat{w}(t), \label{pf:04}
\eequ
where $\hat{w}(t)\in\IR^n$ is an estimate of $w$ obtained from
\bequ
    \dot{\hat{w}}(t)=K\mathcal{Q}(\mathcal{G})\bigl(x(t)-\hat{x}(t)\bigl), \quad \hat{w}(0)=\hat{w}_0, \label{pf:05}
\eequ
with $K=\diag(k) \in \IS_+^{n \times n}$, $k=[k_1,\ldots,k_n]\mT$, being the learning rate matrix,
\bequ
    \mathcal{Q}(\mathcal{G}) = \eye_n-\mathcal{S}(\mathcal{G})\bigl(\eye_n+\mathcal{A}(\mathcal{G})\bigl) \in \IR^{n \times n},  \label{pf:06}
\eequ
$\mathcal{S}(\mathcal{G})=\diag(s)\in \IS_+^{n \times n}$, $s=[(N_1+1)^{-1}, \ldots, (N_n+1)^{-1}]\mT$, with $N_i$, $i=1,\ldots,n$, being the number of neighbors that agent $i$ has,
and $\hat{x}(t)\in\IR^n$ is predicted state vector \cite{ref:005} obtained from
\bequ
    \dot{\hat{x}}(t)=-\mathcal{L}(\mathcal{G})\hat{x}(t)+M\bigl(x(t)-\hat{x}(t)\bigl), \quad \hat{x}(0)=\hat{x}_0, \label{pf:07}
\eequ
with $M=m\eye_n \in \IS_+^{n \times n}$ being the predicted state gain.

In order to express the proposed controller for an agent $i$, we first show that
\bequ
    \mathcal{Q}(\mathcal{G})&=&\mathcal{S}(\mathcal{G})\mathcal{L}(\mathcal{G}). \label{pf:07:X}
\eequ
To see this, we write
\vspace{0.1cm}
\bequ
        \eye_n + \mathcal{A}(\mathcal{G}) &=& \eye_n + \mathcal{A}(\mathcal{G}) \pm \mathcal{D}(\mathcal{G}) \nonumber\\
                                          &=& \eye_n + \mathcal{D}(\mathcal{G}) -\mathcal{L}(\mathcal{G}) \nonumber \\
                                          &=& \mathcal{S}(\mathcal{G})^{-1} - \mathcal{L}(\mathcal{G}), \label{pf:08}
\eequ
where
\bequ
    \mathcal{S}(\mathcal{G})&=&[\eye_n + \mathcal{D}(\mathcal{G})]^{-1}, \label{pf:09}
\eequ
by the definition of $\mathcal{S}(\mathcal{G})$.
Then, using (\ref{pf:08}) in (\ref{pf:06}), we have
\bequ
    \mathcal{Q}(\mathcal{G}) &=& \eye_n-\mathcal{S}(\mathcal{G})\bigl(\eye_n + \mathcal{A}(\mathcal{G})\bigl) \nonumber\\
                             &=& \eye_n-\mathcal{S}(\mathcal{G})\bigl(\mathcal{S}(\mathcal{G})^{-1}-\mathcal{L}(\mathcal{G})\bigl) \nonumber\\
                             &=& \mathcal{S}(\mathcal{G})\mathcal{L}(\mathcal{G}), \label{pf:10}
\eequ
which shows the equality given by (\ref{pf:07:X}).

Now, considering the control input (\ref{pf:02}) along with (\ref{pf:03}), (\ref{pf:04}), (\ref{pf:05}), and (\ref{pf:07}),
and using (\ref{pf:10}) in (\ref{pf:05}), the form of the proposed controller for an agent $i$ becomes
\vspace{0.1cm}
\bequ
    u_i(t)             &=& -\sum_{i \sim j} \bigl(x_i(t)-x_j(t) \bigl) - \hat{w}_i(t), \label{pf:21} \\
    \dot{\hat{w}}_i(t)       &=& -[K\mathcal{S}(\mathcal{G})]_{ii} \sum_{i \sim j} \bigl(\tilde{x}_i(t)-\tilde{x}_j(t) \bigl), \label{pf:22} \\
    \dot{\hat{x}}_i(t) &=& -\sum_{i \sim j} \bigl(\hat{x}_i(t)-\hat{x}_j(t) \bigl) + m \bigl(x_i(t)-\hat{x}_i(t) \bigl), \label{pf:23} \ \ \ \
\eequ
where $\tilde{x}_i(t) \triangleq x_i(t)-\hat{x}_i(t)$ and $[K\mathcal{S}(\mathcal{G})]_{ii}$ is the $i$-th diagonal element of the diagonal matrix $K\mathcal{S}(\mathcal{G})$.

\vspace{0.1cm}

\subsection{Discussion}
Before giving a formal mathematical proof of the proposed controller (\ref{pf:02}) to show that it solves the consensus problem for the multiagent system given by (\ref{pf:01}), 
we give here a discussion about its architecture. 
Specifically, the proposed controller given by (\ref{pf:02}) is constructed using the standard control input (\ref{pf:03}) and the additional control input (\ref{pf:04}).
As we discussed in the preliminaries section, the standard control input solves the consensus problem when $w=\zero_n$.
Here, the purpose using the additional control input is to suppress the effect of disturbances in order to reach a consensus.

The additional control input is determined from (\ref{pf:05}) and (\ref{pf:07}).
In particular, (\ref{pf:07}) serves as an \textit{ideal consensus model} capturing the dynamics given by (\ref{eqn:04}) if $x(t)-\hat{x}(t) \approx \zero_n$.
That is, since the difference between the aggregated state vector and its estimation determines how close $x(t)$ is to the ideal consensus model,
we use $x(t)-\hat{x}(t)$ in (\ref{pf:05}) as an error signal.
Notice that since the additional control input (\ref{pf:04}) is an integration of (\ref{pf:05}),
then can be readily seen that it minimizes this error signal multiplied from left by $K\mathcal{Q}(\mathcal{G})$, i.e., $K\mathcal{Q}(\mathcal{G})\bigl(x(t)-\hat{x}(t)\bigl)$.
Hence, the form of $\mathcal{Q}(\mathcal{G})$ is not arbitrary and plays an important role.

\vspace{0.1cm}

Consider the \textit{projections} defined by
\vspace{0.1cm}
\bequ
        P_{\mathcal{L}(\mathcal{G})} &\triangleq& {\mathcal{L}(\mathcal{G})} ({\mathcal{L}(\mathcal{G})}\mT {\mathcal{L}(\mathcal{G})})^+ {\mathcal{L}(\mathcal{G})}\mT, \label{proj:01} \\
        P_{\mathcal{L}(\mathcal{G})}^\bot &\triangleq& \eye_n - P_{\mathcal{L}(\mathcal{G})}, \label{proj:02}
\eequ
where $P_{\mathcal{L}(\mathcal{G})}$ is a projection onto the column space of $\mathcal{L}(\mathcal{G})$
and $P_{\mathcal{L}(\mathcal{G})}^\bot$ is a projection onto the null space of $\mathcal{L}(\mathcal{G})$.
Note that $P_{\mathcal{L}(\mathcal{G})}^\bot = {\one_n \one_n\mT}/n$ holds.
Furthermore, every solution to (\ref{eqn:04}) evolves in time to $\mathcal{N}(\mathcal{L}(\mathcal{G}))$.
Therefore, choosing $\mathcal{Q}(\mathcal{G})$ in a form similar to the projection $P_{\mathcal{L}(\mathcal{G})}=\eye_n-{\one_n \one_n\mT}/n$ removes the \textit{bias} from the error signal $x(t)-\hat{x}(t)$.
However, if we choose $\mathcal{Q}(\mathcal{G})$ as $\eye_n-{\one_n \one_n\mT}/n$, then realization of the proposed controller requires agent $i$ to have access to the relative state information of all other agents.
In order to overcome this problem, we choose $\mathcal{Q}(\mathcal{G})$ given by (\ref{pf:06}), which can be viewed as the \textit{localized} version of the projection $\eye_n - ({\one_n \one_n\mT}/n)$.
Specifically, $\mathcal{S}(\mathcal{G})$ and $\eye_n + \mathcal{A}(\mathcal{G})$ are the localized versions of $1/n$ and ${\one_n \one_n\mT}$, respectively.
Therefore, the realization of the proposed controller only requires agent $i$ to have access to its own state and to relative state information with respect to its neighbors. 

\section{Stability Analysis}

Let $\tilde{x}(t)\triangleq x(t)-\hat{x}(t)$ and $\tilde{w}(t)\triangleq \hat{w}(t)-w$ be the aggregated state error and the disturbance error vectors, respectively.
Using (\ref{pf:01}), (\ref{pf:05}), and (\ref{pf:07}), we can write
\vspace{0.1cm}
\bequ
    \dot{\tilde{x}}(t) &=& -\mathcal{L}(\mathcal{G})x(t)-\hat{w}(t)+w+\mathcal{L}(\mathcal{G})\hat{x}(t)-M\tilde{x}(t) \nonumber\\
                       &=& -\mathcal{L}(\mathcal{G})\tilde{x}(t)-M\tilde{x}(t)-\tilde{w}(t) \nonumber\\
                       &=& \tilde{A} \tilde{x}(t)-\tilde{w}(t), \quad \tilde{x}(0)=x_0-\hat{x}_0, \label{sa:01} \\
    \dot{\tilde{w}}(t) &=& K\mathcal{Q}(\mathcal{G})\bigl(x(t)-\hat{x}(t)\bigl) \nonumber\\
                       &=& K\mathcal{Q}(\mathcal{G})\tilde{x}(t), \quad \tilde{w}(0)=\hat{w}_0-w, \label{sa:02}
\eequ
where $\tilde{A} \triangleq -\mathcal{L}(\mathcal{G}) - M$.
By defining $e(t)\triangleq \bigl[\tilde{x}\mT(t), \ \tilde{w}\mT(t) \bigl]\mT$, we can equivalently write (\ref{sa:01}) and (\ref{sa:02}) as
\vspace{0.1cm}
\bequ
    \dot{e}(t) = \begin{bmatrix} \tilde{A} & -\eye_n  \\K\mathcal{Q}(\mathcal{G})  & \zero_{n \times n}  \end{bmatrix} e(t), \ \ \ e(0)=\begin{bmatrix} x_0-\hat{x}_0  \\ \hat{w}_0-w  \end{bmatrix}. \label{sa:03} \
\eequ


\subsection{Supporting Lemmas}
The $2n \times 2n$ system matrix
\vspace{0.1cm}
\bequ
        \tilde{A}_0 \triangleq \begin{bmatrix} \tilde{A} & -\eye_n  \\K\mathcal{Q}(\mathcal{G})  & \zero_{n \times n}  \end{bmatrix}, \label{sa:04}
\eequ
in (\ref{sa:03}) plays an important role in the stability analysis.
We need the following supporting lemmas in order to analyze the properties of this matrix.
\vspace{0.1cm}

\begin{lemma} \label{lem:01}
$\tilde{A} = -\mathcal{L}(\mathcal{G}) - M$ in (\ref{sa:04}) is asymptotically stable.
\end{lemma}

\textbf{Proof}. Since $-\mathcal{L}(\mathcal{G})$ is Lyapunov stable, then there exists a $P \in \IS_+^{n \times n}$ satisfying the Lyapunov equation
\bequ
0 &=& -\mathcal{L}(\mathcal{G})\mT P - P \mathcal{L}(\mathcal{G}) + R \nonumber\\
  &=& -\mathcal{L}(\mathcal{G}) P - P\mathcal{L}(\mathcal{G}) + R, \label{sa:05}
\eequ
for a given $R\in\overline{\IS}_+^{\hspace{0.1em} n \times n}$ \cite{ref:003}.
Adding and subtracting $2mP$ to (\ref{sa:05}) yields
\bequ
    0 &=& -\mathcal{L}(\mathcal{G}) P - P\mathcal{L}(\mathcal{G}) + R \pm 2mP \nonumber\\
      &=& -\bigl[\mathcal{L}(\mathcal{G})+m \eye_n\bigl]P-P\bigl[\mathcal{L}(\mathcal{G})+m \eye_n\bigl]+R+2mP \nonumber\\
      &=& \tilde{A}P + P\tilde{A}+R+2mP, \label{sa:06}
\eequ
where it follows from (\ref{sa:06}) that
\bequ
    0 > \tilde{A}P + P\tilde{A} = \tilde{A}\mT P + P \tilde{A}, \label{sa:07}
\eequ
since $R+2mP>0$.
Now, since there exists a $P \in \IS_+^{n \times n}$ such that (\ref{sa:07}) holds, then it follows that $\tilde{A}$ is asymptotically stable \cite{ref:003}. \hfill $\square$

\begin{lemma}[{\hspace{-0.02cm}\cite{ref:006}}] \label{lem:02}
Let $A$, $B \in \IR^{n \times n}$. Then, there is a nonsingular matrix $C\in\IR^{n \times n}$ such that $A=C B C\mT$ if and only if $A$ and $B$ have the same number of positive, negative, and zero eigenvalues.
\end{lemma}

\begin{lemma} \label{lem:03}
$K\mathcal{Q}(\mathcal{G})$ in (\ref{sa:04}) has $n-1$ positive eigenvalues and a zero eigenvalue.
\end{lemma}

\textbf{Proof}. First note that $K\mathcal{Q}(\mathcal{G})=K\mathcal{S}(\mathcal{G})\mathcal{L}(\mathcal{G})$ from (\ref{pf:07:X}).
Let $\lambda \in \spec\bigl(K\mathcal{S}(\mathcal{G})\mathcal{L}(\mathcal{G})\bigl)$ be such that $K\mathcal{S}(\mathcal{G})\mathcal{L}(\mathcal{G}) \eta = \lambda \eta$,
where $\eta\in\IC^n$ and $\eta \neq 0$.
Now, noting that $\bigl(K\mathcal{S}(\mathcal{G})\bigl)^{1/2}\mathcal{L}(\mathcal{G})\eta=\lambda \bigl(K\mathcal{S}(\mathcal{G})\bigl)^{-1/2} \eta$ or, equivalently,
$C \mathcal{L}(\mathcal{G}) C \bigl(C^{-1}\eta\bigl) = \lambda \bigl(C^{-1}\eta\bigl)$, $C \triangleq (K\mathcal{S}(\mathcal{G})\bigl)^{1/2}$, it follows that
$C \mathcal{L}(\mathcal{G}) C \tilde{\eta} = \lambda \tilde{\eta}$, where $(\cdot)^{1/2}$ denotes the unique positive-definite square root and $\tilde{\eta}\triangleq C^{-1} \eta$.
This shows that $\spec\bigl(K\mathcal{S}(\mathcal{G})\mathcal{L}(\mathcal{G})\bigl)=\spec\bigl(C \mathcal{L}(\mathcal{G}) C\bigl)$.
Now, let $\tilde{\mathcal{L}}(\mathcal{G}) \triangleq C \mathcal{L}(\mathcal{G}) C$ and note that $C=C\mT$.
Then, it follows from Lemma \ref{lem:02} that $\tilde{\mathcal{L}}(\mathcal{G})$ and $\mathcal{L}(\mathcal{G})$ have $n-1$ positive eigenvalues and a zero eigenvalue.
From the equivalence $\spec\bigl(K\mathcal{Q}(\mathcal{G})\bigl)=\spec\bigl(C\mathcal{L}(\mathcal{G})C\bigl)$, the result is immediate. \hfill $\square$

Lemma \ref{lem:01} shows that $\tilde{A}$ in (\ref{sa:04}) is asymptotically stable, and hence, every eigenvalue of $\tilde{A}$ has negative real part.
In addition, Lemma \ref{lem:03} shows that $K\mathcal{Q}(\mathcal{G})$ in (\ref{sa:04}) has $n-1$ positive eigenvalues and a zero eigenvalue.
In order to use these two results to analyze (\ref{sa:04}) (and hence, (\ref{sa:03})), we also need the following lemma.

\begin{lemma}[{\hspace{-0.02cm}\cite{ref:007}}] \label{lem:04}
Suppose $\mathcal{Z}(\lambda)=A \lambda^2 + B \lambda + C$ denotes the quadratic matrix polynomial, where $A\in\IR^{n \times n}$ and $C\in\IR^{n \times n}$, and $A$ is nonsingular.
If $B\in\IR^{n \times n}$ is positive-definite, then $\pi_+(\mathcal{Z})=\pi_-(A)+\pi_-(C)$, $\ \pi_-(\mathcal{Z})=\pi_+(A)+\pi_+(C)$, and $\pi_0(\mathcal{Z})=\pi_0(C)$,
where $\pi_+(\mathcal{Z})+\pi_-(\mathcal{Z})+\pi_0(\mathcal{Z})=2n$.
\end{lemma}

\subsection{Main Result}

The following theorem analyzes the solution $e(t)$ of the error system given by (\ref{sa:03}).

\begin{theorem} \label{thm:01}
The solution $e(t)$ of the error system given by (\ref{sa:03}) is Lyapunov stable for all $e_0 \in \IR^{2n}$ and $t \ge 0$, and $e(t) \rightarrow \epsilon \bigl[\one_n\mT, \ -m\one_n\mT \bigl]\mT$ as $t \rightarrow \infty$,
where $\epsilon$ is a constant in $\IR$.
\end{theorem}

\textbf{Proof}. Differentiating (\ref{sa:01}) with respect to time and using (\ref{sa:02}), we can write
\vspace{0cm}
\bequ
    \ddot{\tilde{x}}(t)=\tilde{A}\dot{x}(t)-K \mathcal{Q}(\mathcal{G})\tilde{x}(t). \label{mr:01}
\eequ
Define $e_{\xi1}(t)\triangleq\tilde{x}(t)$ and $e_{\xi2}(t)\triangleq\dot{\tilde{x}}(t)$.
Then, letting $e_\xi(t) \triangleq \bigl[e_{\xi1}\mT(t), \ e_{\xi2}\mT(t) \bigl]\mT$ yields
\vspace{-0.1cm}
\bequ
    \dot{e}_\xi(t) = \begin{bmatrix} \zero_{n \times n} & \eye_n  \\-K\mathcal{Q}(\mathcal{G})  & \tilde{A}  \end{bmatrix} e_\xi(t). \label{mr:02}
\eequ
Notice that the system matrix of (\ref{mr:02}) is equivalent to the system matrix of (\ref{sa:03}) in that the characteristic equations of both matrices are the same and can be given from (\ref{mr:02}) as
\vspace{-0.1cm}
\bequ
    \mathcal{Z}(\lambda) = \lambda^2 \eye_n + \lambda (-\tilde{A}) + K\mathcal{Q}(\mathcal{G}) = 0. \label{mr:03}
\eequ
In order to apply Lemma \ref{lem:04}, we only need to show the positive-definiteness of $-\tilde{A}$.
To see this, note that $\xi\mT \mathcal{L}(\mathcal{G}) \xi \ge 0$.
Then, we can write $\xi\mT \bigl(-\tilde{A}\bigl) \xi = \xi\mT \bigl(\mathcal{L}(\mathcal{G}) + m\eye_n \bigl) \xi = \xi\mT \mathcal{L}(\mathcal{G}) \xi + m \xi\mT \xi>0$, $\xi \neq 0$.
Therefore, $-\tilde{A}$ is positive-definite.

In Lemma \ref{lem:01}, we showed that every eigenvalue of $\tilde{A}$ has negative real part.
This implies that every eigenvalue of $-\tilde{A}$ has positive real part.
Furthermore, in Lemma \ref{lem:03}, we showed that $K\mathcal{Q}(\mathcal{G})$ has $n-1$ positive eigenvalues and a zero eigenvalue.
Therefore, it now follows from Lemma \ref{lem:04} that $\pi_+(\mathcal{Z})=0$, $\pi_-(\mathcal{Z})=2n-1$, and $\pi_0(\mathcal{Z})=1$, where $\mathcal{Z}$ is given by (\ref{mr:03}).
By the equivalency of the system matrices appearing in (\ref{sa:03}) and (\ref{mr:02}), it follows that the system matrix $\tilde{A}_0$ given by (\ref{sa:04}) is Lyapunov stable,
and hence, the solution $e(t)$ of the system (\ref{sa:03}) is Lyapunov stable for all $e_0 \in \IR^{2n}$ and $t \ge 0$.

To prove $e(t) \rightarrow \epsilon \bigl[\one_n\mT, \ -m\one_n\mT \bigl]\mT$ as $t \rightarrow \infty$, first note that
the equilibrium solution $e(t)=e_\rom{e}$ to (\ref{sa:03}) corresponds to $e_\rom{e}\in \mathcal{N}(\tilde{A}_0)$, and hence,
every point in the null space of $\tilde{A}_0$ is an equilibrium point for (\ref{sa:03}) \cite{ref:004}.
Therefore, $e(t) \rightarrow \mathcal{N}(\tilde{A}_0)$ as $t \rightarrow \infty$.
Since
\bequ
    \begin{bmatrix} \tilde{A} & -\eye_n  \\K\mathcal{Q}(\mathcal{G})  & \zero_{n\times n}  \end{bmatrix} \begin{bmatrix} p_1  \\ p_2  \end{bmatrix} = \begin{bmatrix} \zero_n  \\ \zero_n  \end{bmatrix}, \label{mr:04}
\eequ
holds for $p_1 = \epsilon \one_n$ and $p_2 = -\epsilon m \one_n$,
we conclude that $\mathcal{N}(\tilde{A}_0)=\span\{\bigl[\one_n\mT, \ -m\one_n\mT \bigl]\mT\}$.
That is, $e(t) \rightarrow \epsilon \bigl[\one_n\mT, \ -m\one_n\mT \bigl]\mT$ as $t \rightarrow \infty$.
Notice that this implies $\tilde{x}(t) \rightarrow \epsilon \one_n$ and $\tilde{w}(t) \rightarrow -\epsilon m \one_n$ as $t \rightarrow \infty$. \hfill $\square$

The following proposition presents the main result of this section.

\begin{proposition} \label{prop:01}
Consider the multiagent system given by (\ref{pf:01}).
Then, the proposed controller (\ref{pf:02})--(\ref{pf:07}) produces consensus, that is,
\bequ
    \lim_{t\rightarrow \infty} x(t) = \frac{\one_n}{n}\Bigl[\one_n\mT x_0 - \one_n\mT \int_0^t \tilde{w}(\sigma)\rom{d}\sigma \Bigl]. \label{mr:10}
\eequ
\end{proposition}

\textbf{Proof}. Consider the projections defined by (\ref{proj:01}) and (\ref{proj:02}).
We can write
\vspace{0cm}
\bequ
x(t)&=&\eye_n x(t)\nonumber \eequ \bequ
    &=& \bigl(P_{\mathcal{L}(\mathcal{G})} + P_{\mathcal{L}(\mathcal{G})}^\bot \bigl) x(t) \nonumber \\
    &=& x_1(t) + x_2(t),
\eequ
where $x_1(t) \triangleq P_{\mathcal{L}(\mathcal{G})} x(t)$ and $x_2(t) \triangleq P_{\mathcal{L}(\mathcal{G})}^\bot x(t)$.
Therefore,
\bequ
\lim_{t\rightarrow \infty} x(t) = \lim_{t\rightarrow \infty} x_1(t)+  \lim_{t\rightarrow \infty} x_2(t), \label{mr:11}
\eequ
holds.

First, we show that
\bequ
    \lim_{t\rightarrow \infty} x_1(t)=\zero_n. \label{mr:12}
\eequ
Differentiating $x_1(t)$ with respect to time, we have
\bequ
    \dot{x}_1(t) &=& P_{\mathcal{L}(\mathcal{G})} \dot{x}(t) \nonumber\\
                 &=& -P_{\mathcal{L}(\mathcal{G})} \mathcal{L}(\mathcal{G})x(t) - P_{\mathcal{L}(\mathcal{G})}\tilde{w}(t) \nonumber\\
                 &=& - \mathcal{L}(\mathcal{G}) P_{\mathcal{L}(\mathcal{G})} x(t) - P_{\mathcal{L}(\mathcal{G})} \tilde{w}(t) \nonumber\\
                 &=& - \mathcal{L}(\mathcal{G})x_1(t)-P_{\mathcal{L}(\mathcal{G})} \tilde{w}(t). \label{mr:14}
\eequ
Recall from Theorem \ref{thm:01} that $\lim_{t\rightarrow \infty} \tilde{w}(t)=-\epsilon m \one_n$,
and hence, $\lim_{t\rightarrow \infty} P_{\mathcal{L}(\mathcal{G})}\tilde{w}(t)=\zero_n$, since $P_{\mathcal{L}(\mathcal{G})}\one_n = \zero_n$.
Therefore, it follows from $\lim_{t\rightarrow \infty} P_{\mathcal{L}(\mathcal{G})}\tilde{w}(t)=\zero_n$ and $P_{\mathcal{L}(\mathcal{G})}^\bot x_1(t)=0$ that $x_1(t) \rightarrow \zero_n$ as $t \rightarrow \infty$.

Next, we show that
\bequ
    \lim_{t\rightarrow \infty} x_2(t)=\frac{\one_n}{n}\Bigl[\one_n\mT x_0 - \one_n\mT \int_0^t \tilde{w}(\sigma)\rom{d}\sigma \Bigl]. \label{mr:15}
\eequ
The solution $x_2(t)$ can be written as
\bequ
    x_2(t) &=& P_{\mathcal{L}(\mathcal{G})}^\bot x(t) \nonumber\\
           &=& \frac{\one_n \one_n\mT}{n} x(t) \nonumber\\
           &=& \frac{\one_n \one_n\mT}{n} \rom{e}^{-\mathcal{L}(\mathcal{G})t}x_0 \nonumber\\
            && - \frac{\one_n \one_n\mT}{n} \int_0^t \rom{e}^{-\mathcal{L}(\mathcal{G})(t-\sigma)} \tilde{w}(\sigma)\rom{d}\sigma \nonumber\\
           &=& \frac{\one_n \one_n\mT}{n} x_o - \int_0^t \frac{\one_n \one_n\mT}{n} \tilde{w}(\sigma) \rom{d} \sigma,
\eequ
which gives (\ref{mr:15}), where we used the fact $\one_n\mT \rom{e}^{-\mathcal{L}(\mathcal{G})} = \one_n\mT$.

Finally, it follows from (\ref{mr:11}) along with (\ref{mr:12}) and (\ref{mr:15}) that the result is immediate. \hfill $\square$

\subsection{Boundedness of the Control Signal}
The next proposition shows that the proposed controller given by (\ref{pf:02}) along with (\ref{pf:03}) and (\ref{pf:04}) produces a bounded signal to the multiagent system given by (\ref{pf:01}).

\begin{proposition} \label{prop:02}
The control signal $u(t)$ in (\ref{pf:02}) satisfies $\norm{u(t)}_2 \le u^*$ for all $t \ge 0$, where $u^*\in \IR_+$.
\end{proposition}

\textbf{Proof}. The multiagent system given by (\ref{pf:01}) subject to the controller given by (\ref{pf:02}) along with (\ref{pf:03}) and (\ref{pf:04}) can be written as
\bequ
    \dot{x}(t) &=& -\mathcal{L}(\mathcal{G})x(t)-\tilde{w}(t) \nonumber\\
               &=& A x(t) + d(t), \label{ctrl:01}
\eequ
where $A \triangleq -\mathcal{L}(\mathcal{G})$ and $d(t) \triangleq - \tilde{w}(t)$.
By defining $P_\Sigma \triangleq [p_1, \ldots, p_{n-1}, \epsilon \one_n]$, where $p_i \in \IR^n$, $i=1,\ldots,n-1$,
$A$ can be decomposed as $A=P_\Sigma \Sigma P_\Sigma\mT$, where
\bequ
        \Sigma \triangleq \begin{bmatrix} \Sigma_0 & \zero_{n-1}  \\ \zero_{n-1}\mT  & 0  \end{bmatrix}, \label{ctrl:02}
\eequ
with $\Sigma_0$ being asymptotically stable.
By applying the transformation $z(t)=P_\Sigma\mT x(t)$, (\ref{ctrl:01}) becomes
\bequ
        \dot{z}(t) &=& \Sigma z(t) + P_\Sigma\mT d(t),
\eequ
or, equivalently,
\bequ
    \dot{z}_1(t) &=& \Sigma_0 z_1(t) + d_1(t), \label{ctrl:03} \\
    \dot{z}_2(t) &=& d_2(t), \label{ctrl:04}
\eequ
where $z_1(t)\in\IR^{n-1}$, $z_2(t)\in\IR$, and $P_\Sigma\mT d(t) = [d_1\mT(t), \ d_2(t)]\mT$ with $d_1(t)\in\IR^{n-1}$ and $d_2(t)\in\IR$.
Notice that the solution $z_1(t)$ to (\ref{ctrl:03}) is bounded, since $\Sigma_0$ is asymptotically stable and $d(t)=-\tilde{w}(t)$ is bounded from Theorem \ref{thm:01}.

Next, the control signal given by (\ref{pf:02}) can be written as
\bequ
    u(t) &=& - \mathcal{L}(\mathcal{G})x(t)-\hat{w}(t) \nonumber\\
         &=& A x(t) - \hat{w}(t) \nonumber\\
         &=& P_\Sigma \Sigma P_\Sigma\mT x(t) - \hat{w}(t) \nonumber\\
         &=& P_\Sigma \Sigma z(t)-\hat{w}(t) \nonumber\\
         &=& P_\Sigma \begin{bmatrix} \Sigma_0 z_1(t)  \\ 0  \end{bmatrix} - \hat{w}(t), \label{ctrl:05}
\eequ
where the first term in (\ref{ctrl:05}) is bounded since the solution $z_1(t)$ to (\ref{ctrl:03}) is bounded and also the second term in (\ref{ctrl:05}) is bounded since $\tilde{w}(t)$ is bounded and
$\tilde{w}(t)=\hat{w}(t)-w$.
Therefore, there exists a $u^*\in \IR_+$ such that $\norm{u(t)}_2 \le u^*$ for all $t \ge 0$. \hfill $\square$

\subsection{Formation Problem}
Similar to Section II.C, the proposed controller for the formation problem is given by
\bequ
        u(t)=u_\rom{s}(t) + u_\rom{a}(t) + u_\rom{f}(t), \label{fp:01}
\eequ
where $u_\rom{s}(t)$ and $u_\rom{a}(t)$ satisfy (\ref{pf:03}) and (\ref{pf:04}), respectively, and the formation control signal $u_\rom{f}(t)$ satisfies
\bequ
        u_\rom{f}(t) = -\mathcal{L}(\mathcal{G})\zeta. \label{fp:02}
\eequ
Furthermore, in this case, we use
\bequ
    \dot{\hat{x}}(t)&=&-\mathcal{L}(\mathcal{G})\hat{x}(t)+M\bigl(x(t)-\hat{x}(t)\bigl)+u_\rom{f}(t), \nonumber\\ && \hspace{3.56cm} \quad \hat{x}(0)=\hat{x}_0,  \label{fp:03}
\eequ
instead of (\ref{pf:07}).

Note that the aggregated state error and the disturbance error vectors given by (\ref{sa:01}) and (\ref{sa:02}) remain the same for the formation problem.
Therefore, Theorem \ref{thm:01} still holds.
Considering Proposition \ref{prop:01}, it can be easily shown that
\bequ
    \lim_{t\rightarrow \infty} x(t) = \zeta + \frac{\one_n}{n}\Bigl[\one_n\mT x_0 - \one_n\mT \int_0^t \tilde{w}(\sigma)\rom{d}\sigma \Bigl], \label{fp:05}
\eequ
holds by applying the state transformation given by (\ref{eqn:05}).
Finally, it can be shown similar to Proposition \ref{prop:02} that the proposed controller given by (\ref{fp:01}) produces a bounded signal to the multiagent system given by (\ref{pf:01}),
since $u_\rom{f}(t)$ given by (\ref{fp:02}) is bounded.

\section{Convergence to a Constant Point}
In the previous section, we showed that the proposed controller (\ref{pf:02})--(\ref{pf:07}) produces consensus.
However, the agents do not necessarily come to a global agreement on a constant point in space,
since the second term in (\ref{mr:10}) (resp. the third term in (\ref{fp:05})) varies with time.
For the applications when the agents need to come to a global agreement on a constant point in space,
then one needs to modify (\ref{pf:05}),
such that $\tilde{w}(t) \rightarrow \zero_n$ as $t \rightarrow \infty$,
and hence, the second term in (\ref{mr:10}) (resp. the third term in (\ref{fp:05})) converges to a constant,
that is $\lim_{t\rightarrow \infty} x(t)=\one_n \alpha$ holds, where $\alpha$ is a constant in $\IR$.

For this purpose, consider
\bequ
    \dot{\hat{w}}(t)=K\bigl(\mathcal{Q}(\mathcal{G}) + q \eye_n\bigl)\bigl(x(t)-\hat{x}(t)\bigl), \quad \hat{w}(0)=\hat{w}_0, \label{const:01}
\eequ
instead of (\ref{pf:05}), where $q \in \IR_+$ is an additional design parameter (which can be chosen close to zero in practice).
Note that (\ref{const:01}) can be expressed for an agent $i$ as
\bequ
    \dot{\hat{w}}_i(t)       &=& -[K\mathcal{S}(\mathcal{G})]_{ii} \sum_{i \sim j} \bigl(\tilde{x}_i(t)-\tilde{x}_j(t) \bigl) \nonumber\\
                              && +q[K]_{ii} \tilde{x}_i(t), \label{const:02}
\eequ
using the similar arguments given in Section III.A.

Considering (\ref{const:01}), the aggregated state error and the disturbance error vectors given by (\ref{sa:01}) and
\bequ
    \dot{\tilde{w}}(t) &=& K\bigl(\mathcal{Q}(\mathcal{G}) + q \eye_n\bigl)\tilde{x}(t), \quad \tilde{w}(0)=\hat{w}_0-w, \ \ \ \ \label{const:03}
\eequ
respectively.
By defining $e(t)\triangleq \bigl[\tilde{x}\mT(t), \ \tilde{w}\mT(t) \bigl]\mT$, we can equivalently write (\ref{sa:01}) and (\ref{const:03}) as
\bequ
    && \hspace{-0.7cm} \dot{e}(t) = \begin{bmatrix} \tilde{A} & -\eye_n  \\K\bigl(\mathcal{Q}(\mathcal{G}) + q \eye_n\bigl)  & \zero_{n \times n}  \end{bmatrix} e(t), \ e(0)=\begin{bmatrix} x_0-\hat{x}_0  \\ \hat{w}_0-w  \end{bmatrix}, \nonumber\\
    &&\label{const:04}
\eequ
where in this case the $2n \times 2n$ system matrix
\bequ
        \tilde{A}_0 \triangleq \begin{bmatrix} \tilde{A} & -\eye_n  \\K\bigl(\mathcal{Q}(\mathcal{G}) + q \eye_n\bigl)  & \zero_{n \times n}  \end{bmatrix}, \label{const:05}
\eequ
in (\ref{const:04}) plays an important role in the stability analysis.
One can show using similar arguments as in Lemmas \ref{lem:01} and \ref{lem:03} that $-K\bigl(\mathcal{Q}(\mathcal{G}) + q \eye_n\bigl)$ is asymptotically stable.
Hence, by applying Lemma \ref{lem:04} as in the proof of Theorem \ref{thm:01},
it follows that $e(t) \rightarrow \zero_{2n}$ as $t \rightarrow \infty$, since $-K\bigl(\mathcal{Q}(\mathcal{G}) + q \eye_n\bigl)$ and $\tilde{A}$ are asymptotically stable
(and hence, $\tilde{A}_0$ given by (\ref{const:05}) is asymptotically stable).
That is, $\tilde{x}(t) \rightarrow \zero_n$ and $\tilde{w}(t) \rightarrow \zero_n$ as $t \rightarrow \infty$.
Therefore, $\lim_{t\rightarrow \infty} x(t)=\one_n \alpha$ is constant, which is a direct consequence of Proposition \ref{prop:01}.
Note that, it can be shown similar to Proposition \ref{prop:02} that the proposed controller given by (\ref{pf:02}) still produces a bounded signal to the multiagent system given by (\ref{pf:01}).
Finally, the discussion given in Section IV.D holds for the formation problem with (\ref{pf:05}) replaced by (\ref{const:01}).

\section{Time-Varying Disturbances}
This section deals with the case when the disturbances are time-varying.
That is, the dynamics of agents are perturbed by an \textit{unknown} \textit{time-varying} disturbance vector $w(t)=[w_1(t),\ldots,w_n(t)]\mT \in \IR^n$,
where $w_i(t)$, $i=1,\ldots,n$, represents a disturbance effecting the dynamics of $i$-th agent\footnote{For the case when there is no disturbance effect on the dynamics of an agent $i$, $w_i(t)=0$}.
Then, the multiagent system is represented by
\bequ
    \dot{x}(t)=u(t)+w(t), \quad x(0)=x_0. \label{tv:00}
\eequ
Here, we assume that the disturbance vector satisfies $\norm{w(t)}_2 \le w^*$ and $\norm{\dot{w}(t)}_2 \le \dot{w}^*$, where $w^*$, $\dot{w}^* \in \IR_+$.

We consider the same proposed controller as given in Section III.A with (\ref{pf:05}) replaced by
\bequ
    \dot{\hat{w}}(t)=K\Bigl[\mathcal{Q}(\mathcal{G})\bigl(x(t)-\hat{x}(t)\bigl) - \kappa \hat{w} \Bigl], \quad \hat{w}(0)=\hat{w}_0, \label{tv:01}
\eequ
where $\kappa\in\IR^+$ is an additional design parameter that is used to add damping to (\ref{tv:01}).
In this case, the aggregated state error and the disturbance error vectors given by (\ref{sa:01}) and
\bequ
    \dot{\tilde{w}}(t) &=& K\Bigl[\mathcal{Q}(\mathcal{G})\tilde{x}(t) - \kappa \hat{w} \Bigl] - \dot{w}, \ \ \tilde{w}(0)=\hat{w}_0-w, \ \ \ \ \ \label{tv:02}
\eequ
respectively.

\begin{assumption} \label{assume:01}
There exists $R$, $\bar{R} \in \IS_+^{n \times n}$ such that
\bequ
        -R       &\triangleq&  2\tilde{A}+\frac{1}{\mu} \eye_n <0, \label{a:01} \\
        \bar{R}  &\triangleq&  \kappa \eye_n -K^{-1} -\mu \bar{Q}\mT \bar{Q} >0, \label{a:02}
\eequ
holds, where $\bar{Q}\triangleq \mathcal{S}(\mathcal{G})\bigl(\eye_n+\mathcal{A}(\mathcal{G})\bigl)$ and $\mu \in \IR^+$ is an arbitrary constant.
\end{assumption}

\begin{remark} \label{remark:01}
In Assumption \ref{assume:01}, (\ref{a:01}) can equivalently be written as
\bequ
    -R       &=&  2\tilde{A}+\frac{1}{\mu} \eye_n \nonumber\\
             &=& -2\mathcal{L}(\mathcal{G})-2m \eye_n +\frac{1}{\mu} \eye_n \nonumber\\
             &=& -2\mathcal{L}(\mathcal{G})-(2m-\frac{1}{\mu})\eye_n <0, \label{a:03}
\eequ
and hence, the sufficient condition that (\ref{a:01}) (resp. (\ref{a:03})) holds can be given by
\bequ
        m > \frac{1}{2\mu}. \label{a:04}
\eequ
In order to satisfy (\ref{a:02}) for a given $\kappa$ design parameter, one can choose $\lambda_\rom{min}(K)$ to be large and $\mu$ to be small.
Notice that a small $\mu$ leads to a large $m$ by (\ref{a:04}).
Therefore, there always exists a set of design parameters $m$, $\kappa$, and $K$ such that Assumption \ref{assume:01} is satisfied.
\end{remark}

\begin{theorem} \label{thm:02}
Consider the multiagent system given by (\ref{tv:00}).
Let the controller be given by (\ref{pf:02}), (\ref{pf:03}), (\ref{pf:04}), (\ref{pf:07}), and (\ref{tv:01}) subject to Assumption \ref{assume:01}.
Then, the closed-loop error signals given by (\ref{sa:01}) and (\ref{tv:02}) are uniformly bounded for all $(\tilde{x}(0),\tilde{w}(0)) \in \mathcal{D}_\alpha$,
where $\mathcal{D}_\alpha$ is a compact positively invariant set, with ultimate bound $\norm{e(t)}_2 < \varepsilon$, $t \ge T$, where
\bequ
        \varepsilon &>& \bigl[\nu_x^2 + \lambda_\rom{max}(K^{-1}) \nu_w^2 \bigl]^{1/2}, \label{ub:01} \\
        \nu_x &\triangleq& \bigl[c/\lambda_\rom{min}(R) \bigl]^{1/2}, \label{ub:02} \\
        \nu_w &\triangleq& \bigl[c/\lambda_\rom{min}(\bar{R}) \bigl]^{1/2}, \label{ub:03}
\eequ
where $c \triangleq \kappa w^{*2}+\lambda_\rom{max}(K^{-1})\dot{w}^{*2}$.
\end{theorem}

\textbf{Proof}. Consider the Lyapunov-like function candidate
\bequ
        V(\tilde{x},\tilde{w}) &=& \tilde{x}\mT \tilde{x} + \tilde{w}\mT K^{-1} \tilde{w}. \label{oo:01}
\eequ
Note that (\ref{oo:01}) satisfies $\hat{\alpha}(\norm{e}_2)\le V(e) \le \hat{\beta}(\norm{e}_2)$, where $\hat{\alpha}(\norm{e}_2)=\hat{\beta}(\norm{e}_2)=\norm{e}_2^2$
with $\norm{e}_2^2 = \tilde{x}\mT \tilde{x} + \tilde{w}\mT K^{-1} \tilde{w}$.
Furthermore, note that $\hat{\alpha}(\cdot)$ and $\hat{\beta}(\cdot)$ are class $\mathcal{K}_\infty$ functions.
Differentiating (\ref{oo:01}) along the closed-loop system trajectories of (\ref{sa:01}) and (\ref{tv:02}) yields
\bequ
    \dot{V}(\cdot) &=& 2\tilde{x}\mT(t)\tilde{A}\tilde{x}(t)-2\tilde{x}\mT(t) \tilde{w}(t)+2\tilde{w}\mT(t)\mathcal{Q}(\mathcal{G})\tilde{x}(t) \nonumber\\
                    && -2\kappa \tilde{w}\mT(t)\hat{w}(t)-2\tilde{w}\mT(t)K^{-1}\dot{w}(t) \nonumber\\
                   &=& 2\tilde{x}\mT(t)\tilde{A}\tilde{x}(t)-2\tilde{w}\mT(t)\bar{Q}\tilde{x}(t)-2\kappa \tilde{w}\mT(t)\tilde{w}(t) \nonumber\\
                    && -2\kappa\tilde{w}\mT(t) w(t) -2\tilde{w}\mT(t) K^{-1} \dot{w}(t). \label{oo:02}
\eequ
Using the inequalities
\bequ
    |-2\tilde{w}\mT\bar{Q}\tilde{x}| &\le& \mu \tilde{w}\mT\bar{Q}\mT\bar{Q}\tilde{w}+\frac{1}{\mu}\tilde{x}\mT\tilde{x}, \\
    |-2\kappa \tilde{w}\mT w| &\le& \kappa \tilde{w}\mT \tilde{w} + \kappa w\mT w,\\
    |-2\tilde{w}\mT K^{-1} \dot{w}| &\le& \tilde{w}\mT K^{-1} \tilde{w} + \dot{w}\mT K^{-1} \dot{w},
\eequ
in (\ref{oo:02}), it follows that
\bequ
    \dot{V}(\cdot) &\le& \tilde{x}\mT(t)\bigl[2\tilde{A}+\frac{1}{\mu}\eye_n \bigl]\tilde{x}(t) - \tilde{w}\mT(t)\bigl[-\mu \bar{Q}\mT \bar{Q}-\kappa \eye_n \nonumber\\
                    && + K^{-1} +2 \kappa \eye_n \bigl]\tilde{w}(t) + \kappa w\mT(t)w(t)\nonumber\\
                    && + \dot{w}\mT(t)K^{-1} \dot{w}(t) \nonumber\\
                   &=& -\tilde{x}\mT(t)R \tilde{x} - \tilde{w}\mT(t) \bar{R} \tilde{w}(t) + \kappa w\mT(t) w(t) \nonumber\\
                    && + \dot{w}\mT(t)K^{-1} \dot{w}(t) \nonumber\\
                   &\le& -\lambda_\rom{min}(R)\norm{\tilde{x}}_2^2 - \lambda_\rom{min}(\bar{R})\norm{\tilde{w}}_2^2 + c. \label{oo:03}
\eequ
Now, for $\norm{\tilde{x}(t)}_2 \ge \nu_x$ or $\norm{\tilde{w}(t)}_2 \ge \nu_w$, it follows that $\dot{V}(\cdot)\le 0$ for all $(\tilde{x}(t),\tilde{w}(t)) \in \mathcal{D}_e \setminus \mathcal{D}_r$, where
\bequ
        \mathcal{D}_e &\triangleq& \{ (\tilde{x},\tilde{w})\in\IR^n \times \IR^n: x \in \IR^n \}, \\
        \mathcal{D}_r &\triangleq& \{ (\tilde{x},\tilde{w})\in\IR^n \times \IR^n: \norm{\tilde{x}(t)}_2 \ge \nu_x \ \nonumber\\ && \hspace{2.6cm} \mathrm{or} \ \norm{\tilde{w}(t)}_2 \ge \nu_w \}.
\eequ
Finally, define
\vspace{0cm}
\bequ
        \mathcal{D}_\alpha &\triangleq& \{ (\tilde{x},\tilde{w})\in\IR^n \times \IR^n: V(\tilde{x},\tilde{w})\le\alpha \},
\eequ
where $\alpha$ is the maximum value such that $\mathcal{D}_\alpha \subseteq \mathcal{D}_e$, and define
\bequ
        \mathcal{D}_\beta &\triangleq& \{ (\tilde{x},\tilde{w})\in\IR^n \times \IR^n: V(\tilde{x},\tilde{w})\le\beta \},
\eequ
where $\beta > \hat{\beta}(\xi)=\xi^2 = \nu_x^2 + \lambda_\rom{max}(K^{-1}) \nu_w^2$.
To show ultimate boundedness of the closed-loop system (\ref{sa:01}) and (\ref{tv:02}), note that $\mathcal{D}_\beta \subset \mathcal{D}_\alpha$.
Now, since $\dot{V}(\cdot)\le 0$, $t \ge T$, for all $(\tilde{x},\tilde{w})\in \mathcal{D}_e \setminus \mathcal{D}_r$ and $\mathcal{D}_r \subset \mathcal{D}_\alpha$,
it follows that $\mathcal{D}_\alpha$ is positively invariant.
Hence, if $(\tilde{x}(0),\tilde{w}(0))\in \mathcal{D}_\alpha$, then it follows from Corollary 4.4 of \cite{ref:004} that the solution $(\tilde{x}(t),\tilde{w}(t))$ to (\ref{sa:01}) and (\ref{tv:02})
is ultimately bounded with respect to $(\tilde{x},\tilde{w})$ with ultimate bound $\hat{\alpha}^{-1}(\beta)=\beta^{1/2}$, which yields (\ref{ub:01}). \hfill $\square$

\begin{remark} \label{remark:02}
The ultimate bound given by (\ref{ub:01}) can be made small by making $c$ small.
This leads to a small $\kappa$ and a large $\lambda_\rom{min}(K)$.
Recall also that Remark \ref{remark:01} suggests to choose a large $\lambda_\rom{min}(K)$.
\end{remark}

Since this section deals with the disturbances that are assumed to be time-varying and persistent, we can only achieve \textit{approximate} consensus.
In particular, as discussed in Section III.B, since (\ref{pf:07}) serves as an \textit{ideal consensus model} capturing the dynamics given by (\ref{eqn:04}) when $\tilde{x}(t) \approx \zero_n$,
this suggests to make ultimate bound as small as possible by judiciously choosing the design parameters as outlined in Remark \ref{remark:02}.
Furthermore, it can be shown similar to Proposition \ref{prop:02} that the proposed controller produces a bounded signal to the multiagent system given by (\ref{tv:00}).
Finally, the discussion given in Section IV.D holds for the formation problem with (\ref{pf:05}) replaced by (\ref{tv:01}).

\section{Illustrative Examples}
We now present three numerical examples to demonstrate the efficacy of the proposed controller for consensus and formation problems under constant and time-varying disturbances.
Specifically, we consider a cycle graph with six agents subject to the initial conditions $x(0)=[-0.4, \ -0.2, \ 0.0, \ 0.4, \ 0.6, \ 0.8]\mT$.
All initial conditions associated with the proposed controller are set to zero.

\begin{figure}[b!]\center \epsfig{file=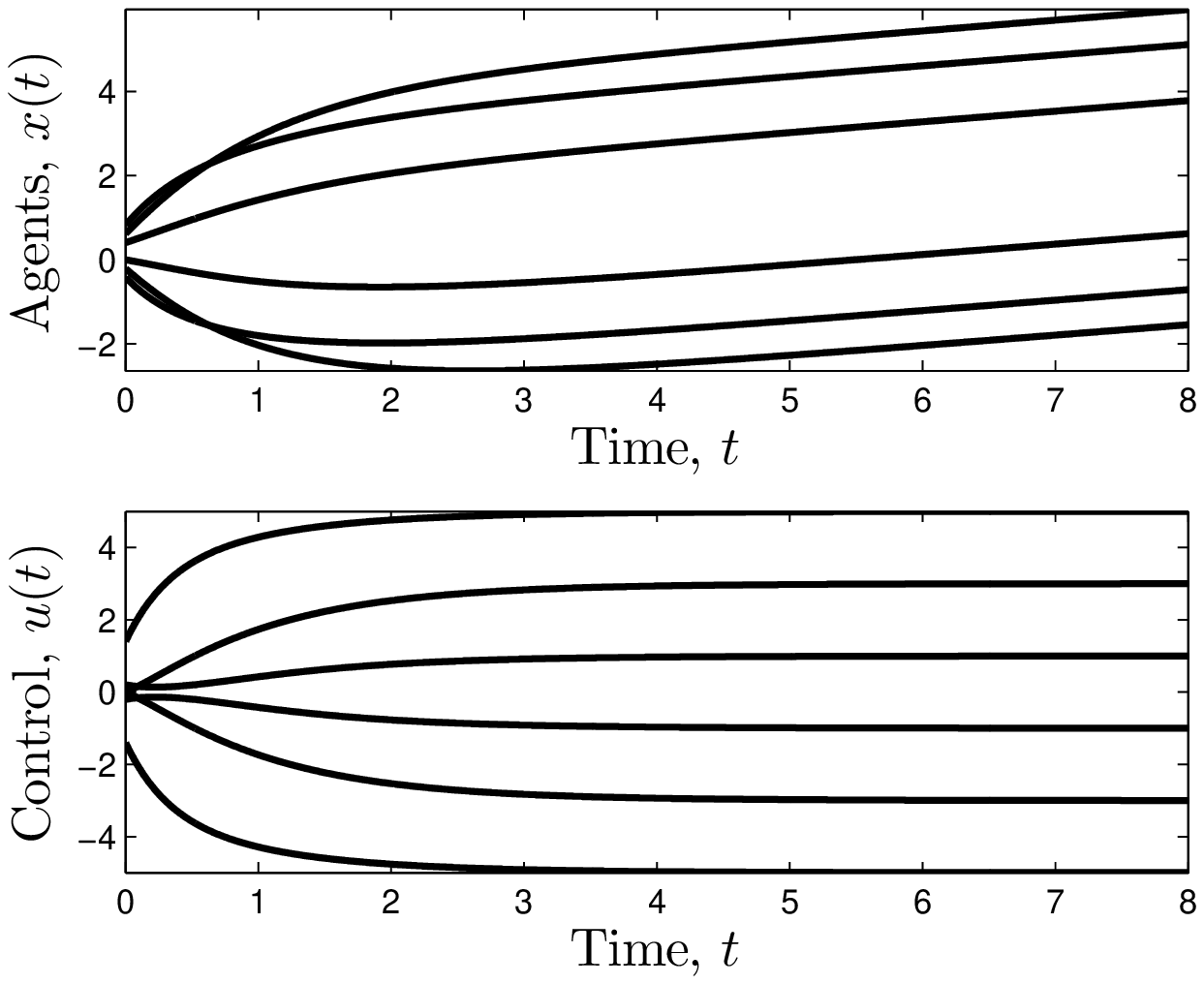, scale=0.5}
\caption{Responses for $x(t)$ and $u(t)$ on the cycle graph with six agents without the additional control input ($u_\rom{a}(t)\equiv0$) under constant disturbances for consensus problem (Example 1).}
\label{fig:1} \end{figure}

\begin{figure}[h!]\center \epsfig{file=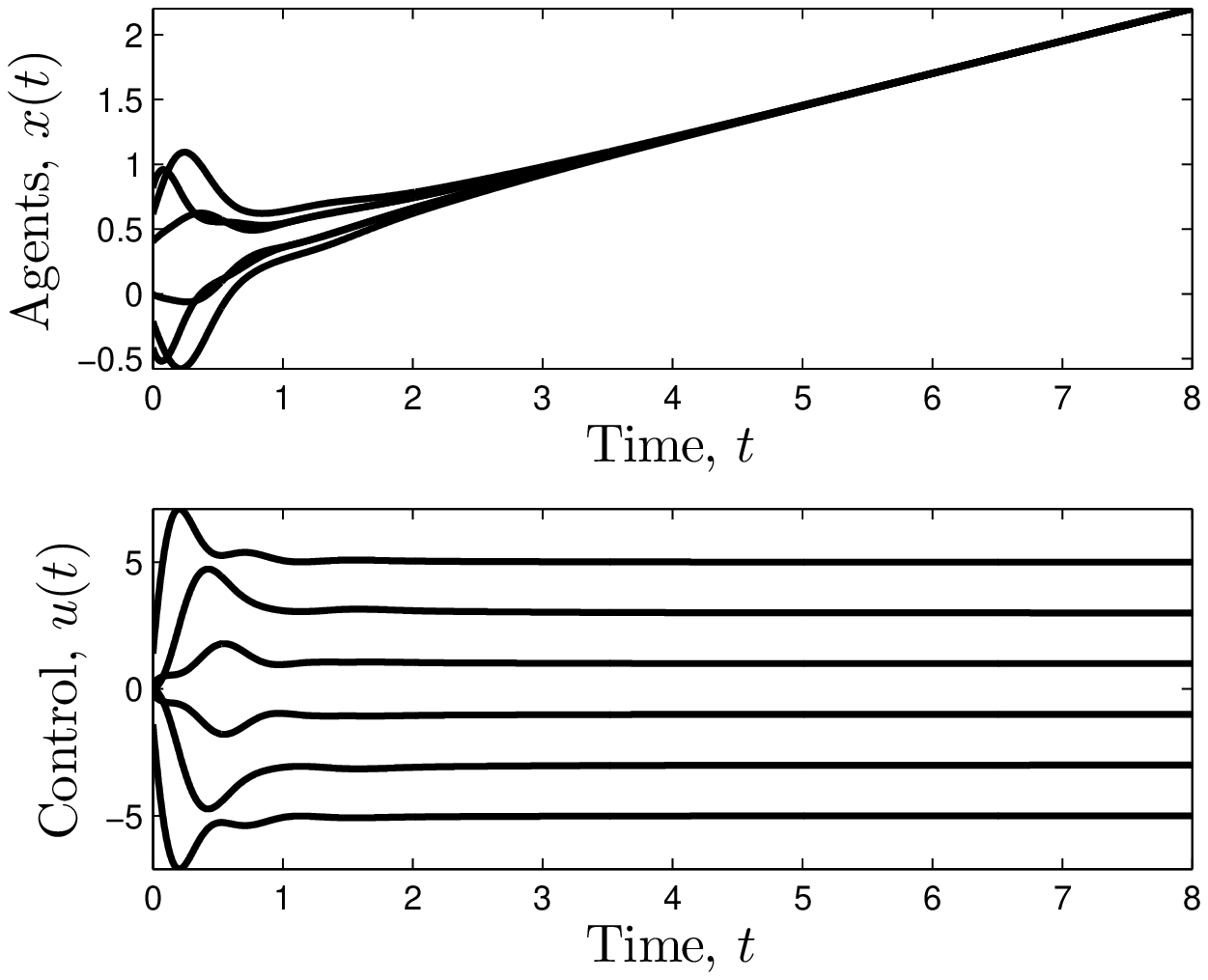, scale=0.5}
\caption{Responses for $x(t)$ and $u(t)$ on the cycle graph with six agents with the additional control input under constant disturbances for consensus problem (Example 1).}
\label{fig:2} \end{figure}

\begin{figure}[h!]\center \epsfig{file=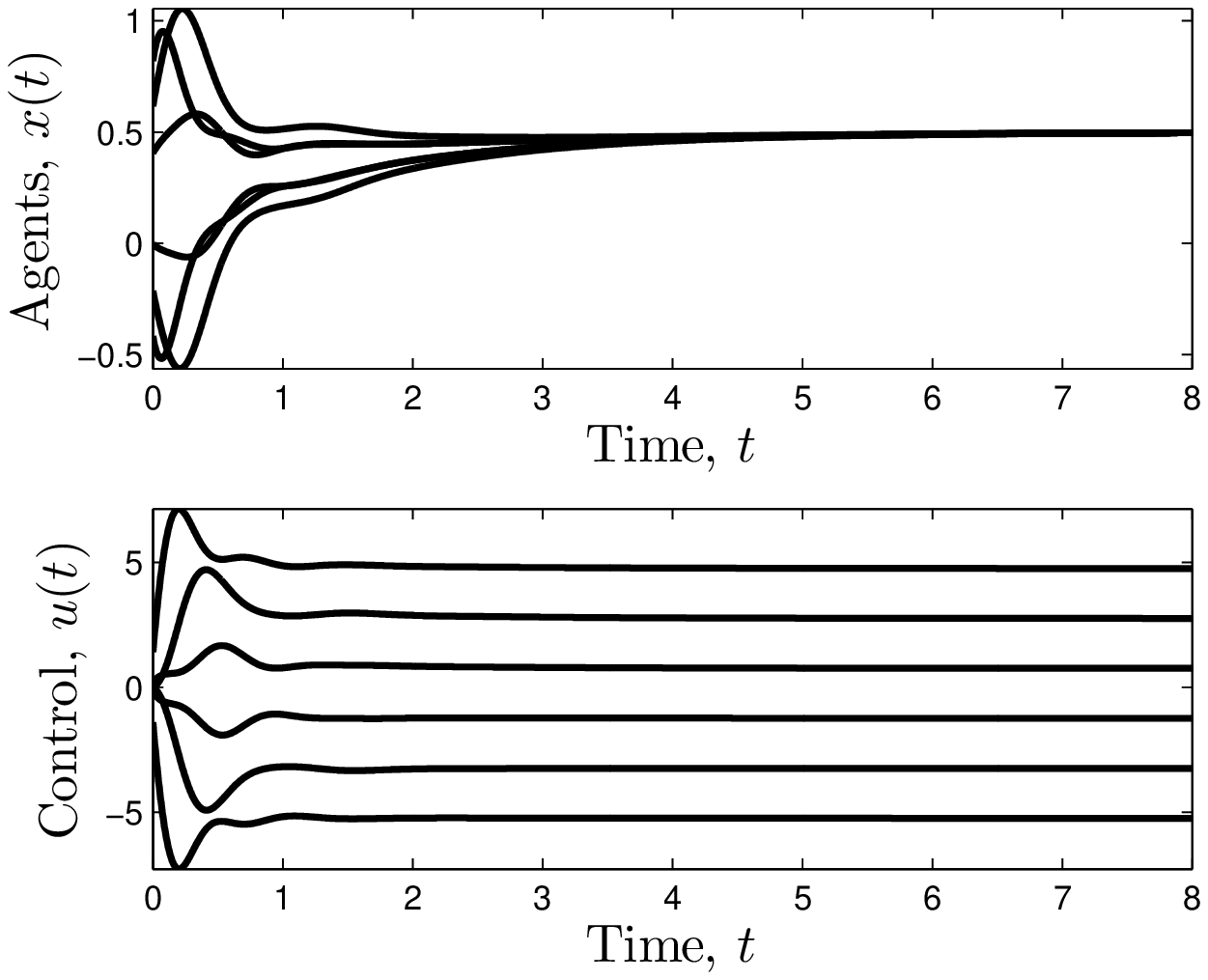, scale=0.5}
\caption{Responses for $x(t)$ and $u(t)$ on the cycle graph with six agents with the additional control input and $q=0.025$ under constant disturbances for consensus problem (Example 1).}
\label{fig:3} \end{figure}

\textit{Example 1: Consensus under Constant Disturbances}.
This example illustrates the consensus problem for the multiagent system given by (\ref{pf:01}) with $w=[-4.75, \ -2.75, \ -0.75, \ 1.25, \ 3.25, \ 5.25]\mT$.
In particular, we use the controller presented in Section III with $K=100 \eye_6$ and $M=5 \eye_6$.
Figure \ref{fig:1} shows the agent positions and control histories when the controller given by (\ref{pf:02}) is applied without the additional control input ($u_\rom{a}(t)\equiv0$).
Figure \ref{fig:2} shows the same histories when both the standard control input (\ref{pf:03}) and the additional control input (\ref{pf:04}) are applied.
The latter figure verifies the presented theory in Section III.
In order to converge to a constant point in space, we employed (\ref{const:01}) with $q=0.025$ (instead of (\ref{pf:05})) in Figure \ref{fig:3},
where this figure verifies the presented theory in Section V. \hfill $\triangle$

\begin{figure}[h!]\center \epsfig{file=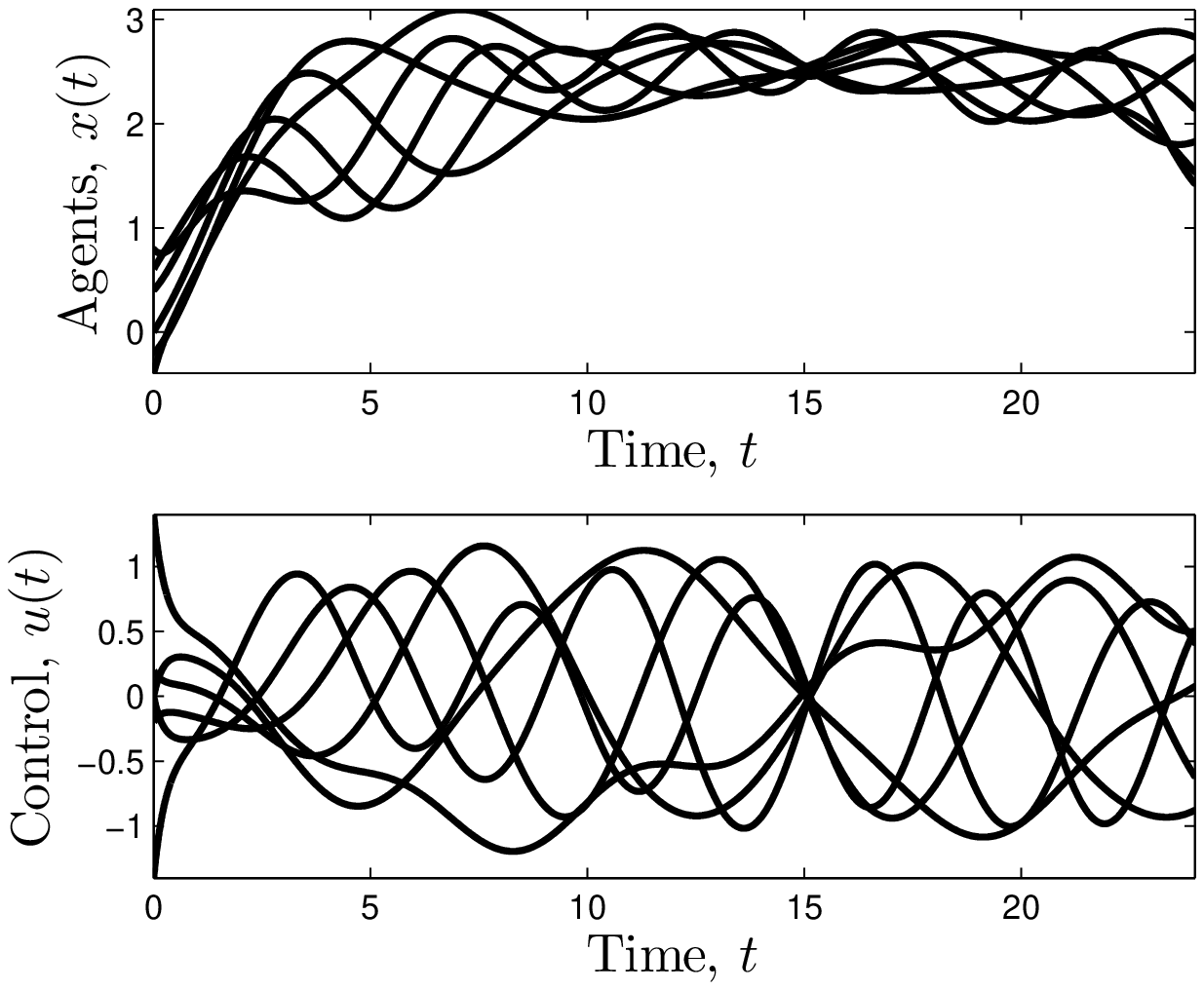, scale=0.5}
\caption{Responses for $x(t)$ and $u(t)$ on the cycle graph with six agents without the additional control input ($u_\rom{a}(t)\equiv0$) under time-varying disturbances for consensus problem (Example 2).}
\label{fig:4} \end{figure}

\begin{figure}[h!]\center \epsfig{file=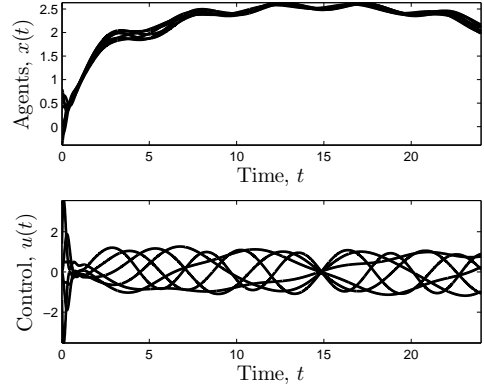, scale=0.5}
\caption{Responses for $x(t)$ and $u(t)$ on the cycle graph with six agents with the additional control input and $\kappa=0.0025$ under time-varying disturbances for consensus problem (Example 2).}
\label{fig:5} \end{figure}

\textit{Example 2: Consensus under Time-Varying Disturbances}.
Next, we illustrate the consensus problem for the multiagent system given by (\ref{tv:00}) with $w(t)=[\rom{sin}(0.2t+10^o), \ \rom{sin}(0.4t+20^o), \rom{sin}(0.6t+30^o), \ \rom{sin}(0.8t+40^o), \ \rom{sin}(1.0t+50^o), \ \rom{sin}(1.2t+60^o)]\mT$.
Here we use the controller presented in Section VI with $K=100 \eye_6$, $M=5 \eye_6$, and $\kappa=0.0025$.
Figure \ref{fig:4} shows the agent positions and control histories when the controller given by (\ref{pf:02}) is applied without the additional control input ($u_\rom{a}(t)\equiv0$).
Figure \ref{fig:5} shows the same histories when both the standard control input (\ref{pf:03}) and the additional control input (\ref{pf:04}) are applied.
The latter figure verifies the presented theory in Section VI.

\hfill $\triangle$

\begin{figure}[h!]\center \epsfig{file=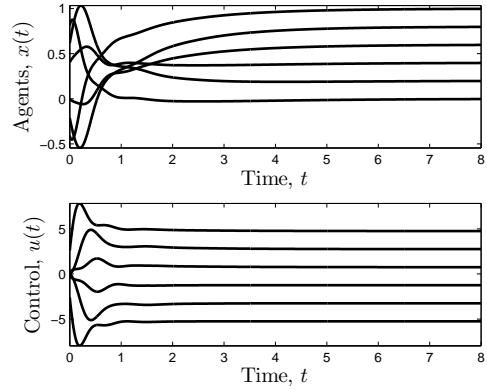, scale=0.5}
\caption{Responses for $x(t)$ and $u(t)$ on the cycle graph with six agents with the additional control input and $q=0.025$ under constant disturbances for formation problem (Example 3).}
\label{fig:6} \end{figure}

\textit{Example 3: Formation under Constant Disturbances}.
Finally, we illustrate the formation problem for the multiagent system given by (\ref{pf:01}) with the same disturbance vector given in Example 1.
We choose $\zeta=[0.0, \ 0.2,\ 0.4,\ 0.6, \ 0.8, \ 1.0]\mT$ for (\ref{fp:02}).
We use the controller presented in Section IV.D with $K=100 \eye_6$ and $M=5 \eye_6$, where (\ref{const:01}) is employed with $q=0.025$ (instead of (\ref{pf:05})).
Figure \ref{fig:6} verifies the presented theory in Section IV.D. \hfill $\triangle$

\section{Conclusions}
To contribute to previous studies of multiagent systems, we have investigated the consensus and formation problems for instances when the dynamics of agents are perturbed by unknown persistent disturbances.
We have shown that the proposed controller suppresses the effect of constant or time-varying disturbances in order to achieve a consensus or a predetermined formation objective.
The realization of the proposed architecture only requires an agent to have access to its own state and to relative state information with respect to its neighbors.
Illustrative examples indicated that the presented theory and its numerical results are compatible.
\\ \\ \\


\end{document}